\documentclass[12pt]{article}
\usepackage{amsfonts,amssymb,amsthm}
\usepackage{pb-diagram}

\newcommand{\vs}{\vspace{5pt}}
\newcommand{\n}{\noindent}
\newcommand{\q}{\quad}
\newcommand{\qq}{\qquad}

\newcommand{\ol}{\overline}

\newcommand{\be}{\beta}
\newcommand{\ga}{\gamma}
\newcommand{\Ga}{\Gamma}
\newcommand{\Th}{\Theta}
\newcommand{\om}{\omega}

\newcommand{\oom}{\ol\om}

\newcommand{\op}{\oplus}
\newcommand{\ot}{\otimes}

\newcommand{\cg}{\mathcal{G}}
\newcommand{\co}{\mathcal{O}}

\newcommand{\g}{\mathfrak{g}}
\newcommand{\fh}{\mathfrak{h}}
\newcommand{\fn}{\mathfrak{n}}
\newcommand{\ft}{\mathfrak{t}}
\newcommand{\cl}{\mathfrak{l}}

\renewcommand{\ge}{\geqslant}

\newcommand{\we}{\wedge}
\newcommand{\We}{\hbox{\small$\bigwedge$}}

\newcommand{\Kur}{\mathrm{Kur}}
\newcommand{\bC}{\mathbb{C}}

\newcommand{\bs}{\backslash}

\renewcommand{\ll}{\langle\!\langle}
\newcommand{\rr}{\rangle\!\rangle}

\newcommand{\half}{\hbox{$\frac12$}}
\newcommand{\bt}{\mathbf{t}}
\newcommand{\bmu}{\hbox{\boldmath$\mu$}}

\newcommand{\bph}{\hbox{\boldmath$\phi$}}
\newcommand{\bPh}{\hbox{\boldmath$\Phi$}}
\newcommand{\dbar}{\overline\partial}
\newcommand{\eps}{\epsilon}
\def\oomega{\overline\omega}

\newcommand{\oper}[2]{\newcommand{#1}{\mathop{\mathrm{#2}}\nolimits} }
\oper{\im}{Im}
\oper{\Aut}{Aut}
\oper{\image}{Im}

\newtheorem{corollary}{Corollary}
\newtheorem{proposition}{Proposition}
\newtheorem{theorem}{Theorem}
\newtheorem{definition}{Definition}

\newtheorem{lemma}{Lemma}
\newcommand{\bproof}{\n{\it Proof: }}
\newcommand{\eproof}{\q QED\vs}
\theoremstyle{remark}
\newtheorem{rem}{Remark}

\begin{document}

\title{Stability of Abelian Complex Structures}

\author{S.~Console$^1$ \and A.~Fino$^1$ \and Y.S.~Poon$^2$}

\footnotetext[1]{Research partially
supported by MIUR (Italy)}
\footnotetext[2]{Partially supported by NSF DMS-0204002.}

\maketitle

\begin{abstract}
Let  $M = \Gamma \backslash G$ be a nilmanifold endowed with an invariant complex structure. We prove that Kuranishi deformations of abelian complex structures  are all  invariant complex structures, generalizing a result  in \cite{MPPS} for 2-step nilmanifolds. We characterize small deformations that remain abelian.  As an application, we observe that at real
dimension six, the deformation process of abelian complex structures is stable
within the class of nilpotent complex structures. We give an example to show that this property
does not hold in higher dimension.

\end{abstract}

\n{AMS Subject Classification: 32G05; 53C15, 53C56, 57S25, 22E25}

\section{Introduction}

Let $M = \Gamma  \backslash G$ be a nilmanifold, i.e. a compact quotient of a
simply-connected nilpotent Lie group $G$ by a uniform discrete subgroup
$\Gamma$. We assume that $M$ has an invariant complex structure
$J$, that is to say that $J$ comes from a complex structure $J$ on the Lie
algebra
$\mathfrak g$ of $G$. An important class of complex structures is given by
the abelian ones \cite{BD,BDM}, which are particular types of nilpotent
complex structures considered in \cite{CFGU}.

If $G$ is 2-step nilpotent or equivalently $M$ is a 2-step nilmanifold,
the deformation of abelian complex structures was studied in \cite{MPPS}, where it is proved that the
Kuranishi process preserves the invariance of the deformed complex
structures, at least for small deformations.

In this notes we show that this result can be generalized to
$k$-step nilmanifolds with abelian complex structures, whatever
$k$ is.

\vs\n{\bf Theorem}\it\q Let $G$ be a simply connected nilpotent
Lie group with co-compact subgroup $\Ga$. Then any abelian
invariant complex structure on $M=\Ga\bs G$ has a locally complete
family of deformations consisting entirely of invariant complex
structures.\rm\vs

For the proof, we construct a family of holomorphic fibrations which
can be derived by the ascending series associated to any nilpotent
Lie algebra (Section~\ref{abeliancxstr}). Then, an inductive argument
shows that the Dolbeault cohomology on $M$ with
coefficients in the structure and holomorphic tangent sheaf can be
computed using
invariant forms and invariant vectors (Lemma~\ref{structure sheaf}
and Theorem~\ref{identification} in Section~\ref{cohom}). So, like in
\cite{MPPS}, one can find harmonic representatives for the  Dolbeault
cohomology on $M$ with
coefficients in the holomorphic tangent sheaf. This allows to prove
that application of Kuranishi's method does not
take one outside the subspace of invariant tensors (Section~\ref{deformation}).

\

It is known that deformation of abelian complex structures is not
stable beginning in real dimension six \cite{MPPS}. We developed
the condition for an infinitesimal deformation to be generated by
a family of abelian complex structures in a coordinate free
manner. Given the results in \cite{MPPS}, it is not surprising to
learn that an element in the first cohomology of the tangent sheaf
is integrable to a family of Abelian complex structures if and
only if it is infinitesimally so. Our computation also
characterizes such elements in terms of the Lie algebra structure
given by the Schouten bracket on the direct sum of the space of
(1,0)-vectors and (0,1)-forms.

\

As an application of the theory developed above, we observe that
at real dimension six, deformation of abelian complex structures
is stable within the class of nilpotent complex structures
(Section~\ref{6-examples}). An example shows that this phenomenon
does not persist in higher dimension (Section~\ref{instability}).

\section{Abelian complex structures}\label{abeliancxstr}

A complex structure $J$ on a Lie algebra $\g$ is called {\it
abelian} if $[JA,JB]=[A,B]$, for all $A,B$ in $\g$ (or,
equivalently,  if the complex space of (1,0)-vectors is an abelian
algebra with respect to Lie bracket) \cite{BD, BDM}. Recall also that
$J$ defines (extending it by left-translation) an invariant almost
complex structure on the group $G$ which is integrable.

The Lie groups or, more generally, the nilmanifolds with
an abelian complex structure are to some extent  dual to complex
parallelizable nilmanifolds: indeed,
  in the complex parallelizable case
  $$
  d \g^{*(1, 0)}  \subset
  \g^{*(2, 0)}\, ,$$
and in the abelian case
$$
d \g^{*(1, 0)} \subset
  \g^{*(1, 1)}\, ,
  $$  where $ \g^{*(p, q)}$ denotes the space of $(p,q)$-forms on
$\g$.

Now assume that the Lie algebra $\g$ is $k$-step nilpotent and set $n
:= \dim _{\bC}\g$.

Let us use like in \cite{CFGU} the ascending series $\{ \g_\ell \}$ with
\[
\begin{array}{l}
\g_0=\{0\}\, , \\
\g_1=\{ X \in \g \mid [X, \g]=0\}\, ,\\
\dots \\
\g_\ell=\{ X \in \g \mid [X, \g]\subset \g_{\ell-1}\}\, ,\\
\dots \\
\g_{k-1}=\{ X \in \g \mid [X, \g]\subset \g_{k-2} \}\, ,\\
\g_k=\g\, .
\end{array}
\]
It is apparent from definition that $\g_\ell/\g_{\ell-1}$ is in
the center of $\g/\g_{\ell-1}$. Since $J$ is abelian $J \g_\ell
\subseteq \g_\ell$. Moreover
\begin{description}
\item[(a)] $\g_\ell/\g_{\ell-1}$ is abelian, and
 \item[(b)]
$\g/\g_{k-1}$ is abelian.
\end{description}
Indeed, since $[\g_{\ell},\g ]\subset \g_{\ell-1}$, we have
$$
[\g_{\ell} / \g_{\ell-1},\g_{\ell} / \g_{\ell-1}] \subset [\g_{\ell}
/ \g_{\ell-1},\g / \g_{\ell-1}] =0.
$$
Moreover, by definition, $\g=\g_{k}$, thus
$\g/\g_{k-1}=\g_k/\g_{k-1}$ is abelian.

Let $\{ \omega _{1}, \bar {\omega }_{1}, \dots ,
\omega _{n},\bar {\omega }_{n} \}$ be a
real basis of $\g ^{\ast }$ satisfying the structure
equations
$$
d\omega _{i} =
\sum \limits _{j,k \leq n} A_{ijk}\,
\omega _{j} \wedge \bar {\omega }_{k}
\quad (1\leq i\leq n) ,
$$ and let $\{ X_{1}, \bar {X}_{1}, \dots ,
X_{n}, \bar {X}_{n} \}$
be the real basis of $\g $ dual to this basis of $1$--forms.
Without loss of generality, we can assume that the basis
$\{ X_{i},\bar {X}_{i}; 1\leq i\leq n \}$ is such that
$\{ X_{n - n_{\ell + 1}}, \bar {X}_{n - n_{\ell +1}}, \dots ,
X_{n} , \bar {X}_{n} \}$
is a real basis of $\g_\ell$, $n_\ell = \dim _{\bC}\g_\ell$.
In fact, proceeding as in the proof of
Theorem~12 in \cite{CFGU}, having chosen a basis
$$\{ X_{1},\bar {X}_{1},\dots ,
X_{n - n_{k-1}},\bar {X}_{n - n_{k-1}} \}$$ of the Lie algebra $\g
/ \g_{k-1}$, we complete it to a basis
$$\{ X_{1},\bar {X}_{1},\dots ,X_{n - n_{k-1}},\bar {X}_{n - n_{k-1}},
X_{n - n_{k-1}+1},\bar {X}_{n - n_{k-1}+1},\dots ,X_{n - n_{k-2}},
\bar {X}_{n - n_{k-2}} \}$$ of the Lie algebra $\g / \g_{k-2}$ and
so on, until we have a basis $$\{ X_{1},\bar {X}_{1},\dots
,X_{n},\bar {X}_{n} \}$$ of the Lie algebra $\g$. Thus,
$$\{ X_{n - n_{\ell + 1}},\bar {X}_{n - n_{\ell +1}},\dots ,
X_{n},\bar {X}_{n} \}$$
is a basis for $\g_\ell,$
  and
$\{ \omega _{i}, \bar {\omega }_{i}; 1 \leq i \leq n-n_\ell \}$
determines the quotient Lie algebra $\g / \g_\ell$. Moreover,

\begin{lemma}\label{dbar closed}
The forms $\bar {\omega}_{1}, \dots ,  \bar {\omega}_{n - n_{\ell}},
\dots , \bar {\omega}_{n}$ are all $\bar \partial$-closed.
\end{lemma}

\section{Cohomology theory}\label{cohom}

As in \cite{CF} \cite{CFGU} \cite{MPPS}, we need to identify
Dolbeault and Lie algebra cohomology, in the spirit of \cite{Nom}.
To this aim we construct a chain of fibrations with tori as fibres
which are associated to the ascending series considered above.

\subsection{Fibrations}\label{fibr}

Let $G$ be a simply connected nilpotent Lie group with Lie algebra
$\g$, $M=\Gamma\backslash G$ be the corresponding nilmanifold,
where $\Gamma$ is a discrete uniform subgroup of $G$. Let $J$ be
an abelian complex structure on $\g$.

Let $\{ \g_\ell \}$ be the ascending series of
Section~\ref{abeliancxstr} associated with $\g$. Then we have the
following exact sequence of Lie algebras.
\[
\begin{array}{lll}
0 \to &\g_{k-1}/\g_{k-2} \to \g/\g_{k-2} \to  &\g/\g_{k-1}\to 0\\
&{\hbox{\rm abelian}}&{\hbox{\rm abelian}}\\
& \dots&\\
0 \to &\g_{\ell}/\g_{\ell-1} \to \g/\g_{\ell-1} \to  &\g/\g_{\ell}\to 0\\
&{\hbox{\rm abelian}}&{}\\
& \dots&\\
0 \to &\g_1 \, \, \longrightarrow \, \, \, \, \g \to  &\g/\g_1\to 0\\
&{\hbox{\rm abelian}}&{}
\end{array}
\]
\begin{rem}\label{split}
We have the following splittings as vector spaces
\[
\begin{array}{ll}
\g/\g_\ell&=\g/\g_{k-1}\oplus \dots \oplus \g_{\ell+2}/\g_{\ell+1}
\oplus \g_{\ell+1}/\g_{\ell}\, , \\
\g/\g_{\ell+1}&=\g/\g_{k-1}\oplus \dots \oplus \g_{\ell+2}/\g_{\ell+1} \, .
\end{array}
\]
In particular, $\g/\g_\ell= \g/\g_{\ell+1} \oplus
\g_{\ell+1}/\g_{\ell}$ as a vector space.
\end{rem}

\medskip

Let $G^\ell$  ($\ell=0,\dots k-1$) and $G^{\ell,\ell-1}$ denote
  the simply connected nilpotent Lie group corresponding to $\g/\g_\ell$
  and $\g_\ell/\g_{\ell-1}$ respectively (note that $G^0=G$). Then we
have the surjective homomorphism (which is actually holomorphic)
$$p_{\ell}: G^{\ell-1} \rightarrow G^\ell$$ with fibre the
abelian Lie group $G^{\ell,\ell-1}$.

Let $\Gamma$ be  the uniform discrete subgroup $\Gamma$ of $G$.
Since $\g_1$ is a rational subalgebra of $\g$, one can see that
$p_1(\Gamma)$ is uniform in $G^1$ \cite{Ra}. So we have the
holomorphic fibration
\[
  \pi_1: \Gamma \backslash G \rightarrow p_1 (\Gamma) \backslash G^1
\]
whose fibre is a torus. Inductively we have holomorphic fibrations
with a torus as fibre
\[ \pi_{\ell- 1}: \Gamma^{\ell-1} \backslash
G^{\ell-1} \rightarrow p_{\ell- 1} ( \Gamma^{\ell-1} ) \backslash
G^\ell
\]
if we set  $\Gamma^{\ell}:= p_{\ell} ( \Gamma^{\ell-1} )$.

For the sake of simplicity we denote $M_\ell:=\Gamma^{\ell}
\backslash G^\ell$. So we have a sequence of  nilmanifolds
$M_\ell$ ($\ell=0,\dots, k-1$) and holomorphic fibrations
$$
\pi_\ell: M_{\ell-1} \to M_{\ell}
$$
whose fibres are tori $T_\ell$, with abelian Lie algebras
$\ft_\ell:=\g_\ell/\g_{\ell-1}$. Note that $M_{k-1}$ is a torus as
well, so that the ``last'' fibration $\pi_{k-1}: M_{k-2} \to
M_{k-1}$ has a torus both as fibre and as base.

\subsection{Lie algebra cohomology}\label{Liealgcohom}

Let $L$ be a Lie group and $\cl$ be its Lie algebra endowed with a
complex structure $J$. Then the complexified Lie algebra has a
type decomposition $\cl_\bC=\cl^{1,0}\op\cl^{0,1}$. These are all
spaces of left-invariant vectors on $L$. The definitions are
extended to invariant $(p,q)$-forms in the standard way. For
instance, $\We^k \cl_{\bC}^{*(0,1)}=\cl^{*(0,k)}$ is the space of
$G$-invariant $(0,k)$-forms.

As in \cite{MPPS} we use  a linear operator $\dbar$ on
$(0,1)$-vectors as follows. For any $(1,0)$-vector $V$ and
$(0,1)$-vector $\bar{U}$, set
\[\dbar_{\bar U}V:=[\bar U,V]^{1,0}.\] We obtain a linear map \[\dbar:
\cl^{1,0}\to\cl^{*(0,1)}\ot\cl^{1,0}.\] In view of Lemma \ref{dbar
closed}, when the complex structure $J$ is abelian, we extend this
definition to a linear map on $\cl^{*(0,k)}\ot\cl^{1,0}$ by
setting
\[\dbar_k(\oom\ot V)=(-1)^k\oom\we\dbar V,\] where $\omega\in
\cl^{*(0,k)}$ and $V\in\cl^{1,0}$. We have a sequence
\[
0\to\cl^{1,0}\to\cl^{*(0,1)}\ot\cl^{1,0}\to\ \cdots\ \to
\cl^{*(0,k-1)}\ot\cl^{1,0} \stackrel{\dbar_{k-1}}\to
\cl^{*(0,k)}\ot\cl^{1,0} \stackrel{\dbar_k}\to\cdots \]  Then we
have (\cite[Lemma 3]{MPPS})

\begin{lemma}\label{seq} The above sequence is a complex, i.e.\
$\dbar_k\circ\dbar_{k-1}=0$ for all $k\ge1$.\end{lemma}

Accordingly we recall the following

\begin{definition}\label{1MPPS} Define $H^k_{\dbar}(\cl^{1,0})$ to be
the $k$-th cohomology  of the above complex; more precisely,
\[H^k_{\dbar}(\cl^{1,0})=\frac{\ker\dbar_k}{\image\dbar_{k-1}}=\frac{\ker\left(\dbar_k:\cl^{*(0,k)}\ot\cl^{1,0}\to
\cl^{*(0,k+1)}\ot\cl^{1,0}\right)}{\dbar_{k-1}\left(\cl^{*(0,k-1)}\ot\cl^{1,0}
\right)}.\] \end{definition}

In the sequel, we will use this cohomology both to the $k$-step
nilpotent Lie algebra $\g$ and the quotients $\g/\g_\ell$. We shall
also need a ``relative" version of Definition \ref{1MPPS} (see the
proof of Theorem~\ref{identification}).

\subsection{Dolbeault cohomology}\label{Dolbeault}

We consider again the fibrations
$$
\pi_\ell: M_{\ell-1} \to M_{\ell}\, , \qquad \ell=0,\dots , k-1\, .
$$
The fibre of $\pi_\ell$ is a torus $T_\ell$, whose abelian Lie
algebra is $\ft_\ell:=\g_\ell/\g_{\ell-1}$.

Recall that the ``last''  fibration $\pi_{k-1}: M_{k-2} \to
M_{k-1}$ has a torus both as fibre and as base.

Our goal is to generalize Lemmata 4 and 5 and Theorem 1 in
\cite{MPPS}. For Lemma 5 and Theorem 1, the idea is to start with the
``last''  fibration which is analogue to the fibration considered in
\cite{MPPS}, and go on inductively. At any step the  basis of the
fibre bundle is not a torus but it has ``good properties" (since it
is the total space of the fibre bundle in the previous step).

For generalizing Lemma 4 in \cite{MPPS}, there is no difficulty,
since all we need is the basis constructed in Section \ref{abeliancxstr}:

\begin{lemma} Let $\co_{M_{\ell-1}}$ and $\Th_{M_{\ell-1}}$ be the
structure sheaf and the tangent
sheaf of ${M_{\ell-1}}$. For $j\ge1$, the direct image sheaves with
respect to ${\pi_\ell}$ are
\[\begin{array}{l}
R^j{\pi_\ell}_*\co_{M_{\ell-1}}=\We^j\ft_\ell^{*(0,1)}\ot\co_{M_{\ell}}=\ft_\ell^{*(0,j)}\ot
\co_{M_{\ell}},\\[5pt]
R^j{\pi_\ell}_*{\pi_\ell}^*\Th_{M_{\ell}}=\We^j\ft_\ell^{*(0,1)}\ot
\Th_{M_{\ell}} =\ft_\ell^{*(0,j)}\ot
\Th_{M_{\ell}}.\end{array}\]\end{lemma}

As a start, let us consider the structure sheaf of a general ${M_{\ell}}$:

\begin{lemma}\label{structure sheaf} Let $\ell=0,\dots , k-1$ and
$\co_{M_{\ell}}$  be the
structure sheaf ${M_{\ell}}$. Then
$$
H^j({M_{\ell}},
\co_{M_{\ell}})=\We^j(\g/\g_\ell)^{*(0,1)}=(\g/\g_\ell)^{*(0,j)}\, .
$$
\end{lemma}

\bproof We prove this lemma by induction, beginning with
$\ell=k-1$ and finishing at $\ell=0$. Note that $M_{k-1}$ is a
torus.

\

 \noindent {\em 1st inductive step:} Consider the ``last"
fibration $\pi_{k-1}: M_{k-2} \to M_{k-1}$ and the corresponding
exact sequence of Lie algebras
\[
\begin{array}{lll}
0 \to &\ft_{k-1}=\g_{k-1}/\g_{k-2} \to \g/\g_{k-2} \to
&\g/\g_{k-1}\to 0\, . \\
&{\hbox{\rm abelian}}&{\hbox{\rm abelian}}
\end{array}
\]
Using  the Leray spectral sequence with respect to the
$\dbar$-operator and the holomorphic projection $\pi_{k-1}$, we
have
\[
E_2^{p,q}=H^p(M_{k-1}, R^q{\pi_{k-1}}_*\co_{M_{k-2}}),\qq
E_{\infty}^{p, q}\Rightarrow H^{p+q}(M_{k-2},\co_{M_{k-2}}).\]
 From the previous lemma, when $q\ge1$, \begin{eqnarray*}
E_2^{p,q} &=&H^p(M_{k-1},\We^q \ft_{k-1}^{*(0,1)}
\ot\co_{M_{k-1}})=\We^q \ft_{k-1}^{*(0,1)} \ot H^p({M_{k-1}},
\co_{M_{k-1}})\\ &=&\We^q \ft_{k-1}^{*(0,1)}\ot\We^p
(\g/\g_{k-1})^{*(0,1)}.\end{eqnarray*} Every element in
$E_2^{p,q}$ is a linear combination of the tensor products of
vertical $(0,q)$-forms and $(0,p)$-forms lifted from the base.
Since these forms are globally defined and the differential $d_2$
is generated by the $\dbar$-operator, we have $d_2=0$. It follows
that the Leray spectral sequence degenerates at the $E_2$-level.
Therefore, since by Remark~\ref{split} $\g/\g_{k-2}= \g/\g_{k-1}
\oplus  \g_{k-1}/\g_{k-2}$ as a vector space,
\begin{eqnarray*}
H^j({M_{k-2}},\co_{M_{k-2}})&=&\bigoplus_{p+q=j}E_2^{p,
q}=\bigoplus_{p+q=j}  \We^q\ft_{k-1}^{*(0,1)} \ot
\We^p(\g/\g_{k-1})^{*(0,1)} \\
&=&\bigoplus_{p+q=j}  \We^q( \g_{k-1}/\g_{k-2})^{*(0,1)} \ot
\We^p(\g/\g_{k-1})^{*(0,1)}\\
&=&\We^j (\g/\g_{k-2})^{*(0,1)}=(\g/\g_{k-2})^{*(0,j)}.
\end{eqnarray*}

\noindent {\em General inductive step:} Consider the fibration
$\pi_{\ell+1}: M_{\ell} \to M_{\ell+1}$ and the corresponding
exact sequence of Lie algebras
\[
\begin{array}{lll}
0 \to &\ft_{\ell+1}=\g_{\ell+1}/\g_{\ell} \to \g/\g_{\ell} \to
&\g/\g_{\ell+1}\to 0\\
&{\hbox{\rm abelian}}&{}
\end{array}
\]
Using again the Leray spectral sequence and the same arguments
$$
E_2^{p,q}=\We^q \ft_{\ell+1}^{*(0,1)} \ot H^p({M_{\ell+1}},
\co_{M_{\ell+1}})\, .
$$
By inductive hypothesis $H^p({M_{\ell+1}},
\co_{M_{\ell+1}})=\We^p(\g/\g_{\ell+1})^{*(0,1)}$, so
\begin{eqnarray*}
H^j({M_{\ell}},\co_{M_{\ell}})&=&\bigoplus_{p+q=j}E_2^{p,
q}=\bigoplus_{p+q=j}  \We^q\ft_{\ell+1}^{*(0,1)} \ot
\We^p(\g/\g_{\ell+1})^{*(0,1)} \\
&=&\bigoplus_{p+q=j}  \We^q( \g_{\ell+1}/\g_{\ell})^{*(0,1)} \ot
\We^p(\g/\g_{\ell+1})^{*(0,1)}\\
&=&\We^j (\g/\g_{\ell})^{*(0,1)}=(\g/\g_{\ell})^{*(0,j)}.
\end{eqnarray*}
\eproof

\begin{theorem}\label{identification} Let $M$ be a $k$-step nilmanifold
with an abelian complex structure.
Let $\ell=0,\dots , k-1$ and  $\Th_{M_{\ell}}$  be the
tangent sheaf ${M_{\ell}}$. Then
$$
H^j({M_{\ell}}, \Th_{M_{\ell}})\cong H^j_{\dbar}((\g/\g_\ell)^{(1,0)})\, .
$$
In particular, there is a natural isomorphism
$H^j(M,\Th_M) \cong H^j_{\dbar}(\g^{1,0})$.\end{theorem}

\bproof Again, we prove this lemma by induction, beginning with
$\ell=k-1$ and finishing at $\ell=0$.

\noindent {\em 1st inductive step:} Consider the ``last" fibration
$\pi_{k-1}: M_{k-2} \to M_{k-1}$ and use the Leray spectral
sequence of $\pi^*_{k-1} \Th_{M_{k-1}}$ because both the base and
the fiber are tori. We are in same setting as in \cite[Lemma 5 and
Theorem 1]{MPPS}. Thus, we have an analogue of \cite[Lemma
5]{MPPS}. i.e.
\[
H^j(M_{k-2}, \pi^*_{k-1}
\Th_{M_{k-1}})=\We^j(\g/\g_{k-2})^{*(0,1)}\ot
(\g/\g_{k-1})^{1,0}\, .
\]
Using the exact sequence
  \[
  0\to (\g_{k-1}/\g_{k-2})^{1,0}\ot
\co_{M_{k-2}} \to \Th_{M_{k-2}} \to \pi_{k-1}^*\Th_{M_{k-1}}\to 0
\]
as in \cite[Theorem 1]{MPPS}, we get
\[
H^j({M_{k-2}}, \Th_{M_{k-2}})\cong
H^j_{\dbar}((\g/\g_{k-2})^{(1,0)}).
\]

\noindent {\em 2nd inductive step:}  Consider the fibration
$\pi_{k-2}: M_{k-3} \to M_{k-2}$. Now the basis is not necessarily
a torus.  The Leray spectral sequence of $\pi^*_{k-2}
\Th_{M_{k-2}}$ yields
\begin{eqnarray*}
E_2^{p,q}&=&\We^q(\g_{k-2}/\g_{k-3})^{*(0,1)}\ot H^p (M_{k-2},
\Th_{M_{k-2}})\\
&=&\We^q(\g_{k-2}/\g_{k-3})^{*(0,1)}\ot
H^p_{\dbar}((\g/\g_{k-2})^{1,0})
\end{eqnarray*}
by the first inductive step.

Elements in $\We^q(\g_{k-2}/\g_{k-3})^{*(0,1)}$ are pulled back to
globally defined
$(0,q)$-forms on $M_{k-3}$. Elements of
$H^p_{\dbar}((\g/\g_{k-2})^{1,0})$ yield globally defined forms  on
$M_{k-3}$ and globally defined sections of $\pi^*_{k-3}
\Th_{M_{k-2}}$. Moreover, $d_2 =0$, since $d_2$ is generated by
$\dbar$. So the spectral sequence degenerates at $E_2$. Thus
\begin{eqnarray}\label{eq:Th}
H^j(M_{k-3}, \pi^*_{k-2} \Th_{M_{k-2}})=\bigoplus_{p+q=j}E_2^{p,
q}=\bigoplus_{p+q=j}  \We^q(\g_{k-2}/\g_{k-3})^{*(0,1)} \ot
H^p_{\dbar}((\g/\g_{k-2})^{1,0}) .
\end{eqnarray}
The latter is the cohomology of the complex
\[
\begin{diagram}
\node{\bigoplus_{p+q=j} } \node{(\g_{k-2}/\g_{k-3})^{*(0,q)}}
\arrow{s,l}{0} \node{ \ot}\node{ \hskip1.8cm
((\g/\g_{k-2})^{*(0,p)}\ot (\g/\g_{k-2})^{1,0} )    } \arrow{s,r}
{\dbar}
\\
\node{\bigoplus_{p+q=j} } \node{ (\g_{k-2}/\g_{k-3})^{*(0,q)}} \node{
\ot}\node{  \hskip2cm    ((\g/\g_{k-2})^{*(0,p+1)}\ot
(\g/\g_{k-2})^{1,0} )    }
\end{diagram}
\]

 Note that
$\g_{k-2}/\g_{k-3}$ is in the center of $\g/\g_{k-3}$. If $V$ is
in $\g/\g_{k-2}$ and $\bar{X}_j$ is in $\g/\g_{k-3}$, $[\bar X_j,
V]_{\g/\g_{k-3}}$ is well defined with respect to the induced Lie
bracket on the quotient algebra $\g/\g_{k-3}$. Recall that
$\{\bar\omega_j: 1\leq j\leq n_{\ell}\}$ forms a basis for
$\g/\g_{\ell}$. Its dual basis is $\{\bar{X}_j: 1\leq j\leq
n_{\ell}\}$. Now we consider a linear map $\dbar'$ on
 as follows.
For $\bar\Omega$ in $(\g/\g_{k-3})^{*(0,j)}$ and $V$ in
$(\g/\g_{k-2})^{1,0}$, define
\begin{equation}\label{eq:dbar'}
\dbar'(\bar \omega \ot V)=\sum_{j=1}^{n-n_{k-3}} (-1)^j \bar
\omega \wedge\bar \omega_j \ot  \pi_{k-3} ([\bar X_j,
V]_{\g/\g_{k-3}}^{1,0}).
\end{equation}
It yields a linear map.
\begin{equation}\label{eq:dbar'-complex}
\dbar': {(\g/\g_{k-3})^{*(0,j)} \ot ( \g/\g_{k-2})^{1,0}}\to
{(\g/\g_{k-3})^{*(0,j+1)} \ot ( \g/\g_{k-2})^{1,0}}
\end{equation}
Now, using the fact that $J$ is an abelian complex structure, one
could verify that $\dbar'\circ\dbar'=0$. Moreover, if $\bar{X}_j$
is dual to an element of $\g_{k-2}/\g_{k-3}$, then $[\bar X_j,
V]_{\g/\g_{k-3}}=0$. Therefore,  the complex with $\dbar'$ and the
complex $0\otimes\dbar$ in the previous paragraph agree.

 We
denote by $H^j_{\dbar}((\g/\g_{k-3}, \g/\g_{k-2})^{1,0})$ the
cohomology of the above $\dbar'$-complex. (It is a sort of
``relative" cohomology.)

Next, we can describe  the $\dbar$-complex on $\g/\g_{k-3}$ using
the splitting
\[
\g/\g_{k-3}=\g/\g_{k-2}\oplus\g_{k-3}/\g_{k-2}=\g/\g_{k-2}\oplus\ft_{k-2}.
\]
It yields
\[
\begin{diagram}
\node{((\g/\g_{k-3})^{*(0,j)}\ot (\g/\g_{k-2})^{1,0})}
\arrow{s,l}{\dbar'}  \arrow{ese,t,3}{\dbar''}\node{\hskip .3cm \oplus
}\node{\hskip .1cm ((\g/\g_{k-3})^{*(0,j)}\ot \ft_{k-2}^{1,0})}
\arrow{s,l}{0}
\\
\node{((\g/\g_{k-3})^{*(0,j+1)}\ot (\g/\g_{k-2})^{1,0})} \node{
\hskip .4cm \oplus}\node{\hskip .1cm ((\g/\g_{k-3})^{*(0,j+1)}\ot
\ft_{k-2}^{1,0})}
\end{diagram}
\]
where $\dbar''$ is the ``$\ft_{k-2}^{1,0}$-component " of $\dbar$ (in
a similar sense as in (\ref{eq:dbar'})).

Now, using the $\dbar'$-complex (\ref{eq:dbar'-complex}) and the
``relative" cohomology $H^j_{\dbar}((\g/\g_{k-3},
\g/\g_{k-2})^{1,0})$, (\ref{eq:Th}) becomes
$$
H^j(M_{k-3}, \pi^*_{k-3} \Th_{M_{k-2}})\cong
H^j_{\dbar}((\g/\g_{k-3}, \g/\g_{k-2})^{1,0})\, .
$$
Consider now the exact sequence
  \[
  0\to (\g_{k-2}/\g_{k-3})^{1,0}\ot
\co_{M_{k-3}} \to \Th_{M_{k-3}} \to \pi_{k-3}^*\Th_{M_{k-2}}\to 0\, ,
\]
or (setting as usual $\ft_{k-2}:=\g_{k-2}/\g_{k-3}$)
  \[
  0\to \ft_{k-2}^{1,0}\ot
\co_{M_{k-3}} \to \Th_{M_{k-3}} \to \pi_{k-3}^*\Th_{M_{k-2}}\to 0\, ,
\]
which induces the long exact sequence
\begin{eqnarray*}
\dots &\to& \ft_{k-2}^{1,0}\ot H^j(M_{k-3},\co_{M_{k-3}})\to
H^j(M_{k-3},\Th_{M_{k-3}}) \to
H^j(M_{k-3},\pi_{k-3}^*\Th_{M_{k-3}}) \\
&\stackrel{\delta_j}\to&\ft_{k-2}^{1,0}\ot
H^{j+1}({M_{k-3}},\co_{M_{k-3}}) \to \dots
\end{eqnarray*}
By the previous results this sequence can be written as
\begin{eqnarray*}
\dots &\to& (\g/\g_{k-3})^{*(0,j)}\ot\ft_{k-2}^{1,0}\to
H^j(M_{k-3},\Th_{M_{k-3}}) \to
H^j_{\dbar}((\g/\g_{k-3}, \g/\g_{k-2})^{1,0})\\
&\stackrel{\delta_j}\to&(\g/\g_{k-3})^{*(0,j+1)}\ot\ft_{k-2}^{1,0}
\to \dots.
\end{eqnarray*}

To compute the coboundary map $\delta_j$, we chase diagram and
find that it is precisely the map $\dbar''$ on the ``relative"
cohomology $H^j_{\dbar}((\g/\g_{k-3}, \g/\g_{k-2})^{1,0})$. Given
the long exact sequence we have
\[H^j(M_{k-3},\Th_{M_{k-3}})\cong\ker
\delta_j\op\frac{(\g/\g_{k-3})^{*(0,j)}\ot\ft_{k-2}^{1,0}}{\delta_{j-1}(H^{j-1}_{\dbar}((\g/\g_{k-3},
\g/\g_{k-2})^{1,0})}.\]\ Let us compare this with
$H^j_{\dbar}((\g/\g_{k-3})^{1,0})$, which is the cohomology of the
complex determined by $\dbar=\dbar'+\dbar''$. Let $[b+t] \in
H^j_{\dbar}((\g/\g_{k-3})^{1,0})$, with $b \in
(\g/\g_{k-3})^{*(0,j)}\ot (\g/\g_{k-2})^{1,0}$ and $t \in
(\g/\g_{k-3})^{*(0,j)}\ot \ft_{k-2}^{1,0}$. Then $b \in \ker \dbar
= \ker \dbar' \cap \ker \dbar''$. Moreover $b +t$ is cohomologous
to $b'+t$ if and only if $b-b'\in \image \dbar=\image \dbar'$.
Thus
$$
\frac{ \ker \dbar' \cap \ker \dbar''}{\image \dbar'}
$$
determines the cohomology class of $b$. This is also the space
$\ker \delta_j$.

Next, observe further that $b+t$ is cohomologous to $b+t'$ if and
only if $t-t'\in \image \dbar=\image \dbar''$, i.e., $t-t'\in
\dbar'' (b'')$. Since $b''$ is mapped by $\dbar$ into
$(\g/\g_{k-3})^{*(0,j)}\ot\ft_{k-2}^{1,0}$, $\dbar' b''=0$, so
$b'' \in \ker \dbar'$. Hence $t-t' \in \dbar'' (\ker
\dbar')=\delta_{j-1} ( H^{j-1}_{\dbar}((\g/\g_{k-3},
\g/\g_{k-2})^{1,0}))$. Therefore
\[
H^j_{\dbar}((\g/\g_{k-3})^{1,0})=\ker \delta_j \oplus
\frac{(\g/\g_{k-3})^{*(0,j)}\ot\ft_{k-2}^{1,0}}{\delta_{j-1} (
H^{j-1}_{\dbar}((\g/\g_{k-3}, \g/\g_{k-2})^{1,0}))}.
\]
It is therefore isomorphic to $H^j(M_{k-3}, \Th_{M_{k-3}})$ as
claimed.

With a shift of indices, it established an inductive process to
complete the proof of the theorem. \eproof

\section{Deformation theory}\label{deformation}

We will show that one can find harmonic representatives in the
Dolbeault cohomology groups. To this goal, we introduce an
invariant Hermitian metric on $M$. We choose such a metric so that
if $\{ X_{1},\bar {X}_{1},\dots ,X_{n},\bar{X}_{n} \}$ is the real
basis of $\g$ constructed in Section~\ref{abeliancxstr},  $\{
Y_{1},J Y_{1},\dots ,Y_{n},JY_{n} \}$ forms a Hermitian frame,
where we set $Y_i:=\frac 1 2 (X_i+\bar{X}_i)$.  We use the
resulting inner product on $\g^{*(0,k)}\ot\g^{1,0}$ to define the
orthogonal complement of $\image \dbar_{k-1}$ in $\ker\dbar_k$.
Denote this space by $\image^\perp\dbar_{k-1}$. Then, with the
same proof as in \cite[Theorem 3]{MPPS}, we find that the harmonic
theory is reduced to finite dimensional linear algebra.

\begin{theorem}\label{harmonic rep} The space $\image^\perp\dbar_{k-1}$ is
a space of harmonic representatives for the Dolbeault cohomology
$H^k(M,\Th_M)$ on the compact complex manifold $M$. In addition,
let $\bmu\in\g^{*(0,k)}\ot\g^{1,0}$. Then $\dbar^*\bmu$ with
respect to the $L_2$-norm on the compact manifold $M$ is equal to
$\dbar^*\bmu$ with respect to the Hermitian inner product on the
finite-dimensional vector spaces
$\g^{*(0,k)}\ot\g^{1,0}$.\end{theorem}

Next one can consider the Schouten-Nijenhuis bracket
\[\{\cdot,\cdot\}:
H^1(X,\Th_X)\times H^1(X,\Th_X)\to H^2(X,\Th_X).\]
Recall that it can be defined as follows. Let $\oom\ot V$ and
$\oom'\ot V'$ be vector-valued (0,1)-forms  representing
elements in $H^1(X,\Th_X)$. Then \[\{\oom\ot V,\oom'\ot V'\}=\oom'\we
L_{V'}\oom\ot V +\oom\we
L_V\oom'\ot V' +\oom\we\oom' \ot [V, V'].\] Using the fact that the
complex structure is abelian, so
$[V, V']=0$ for all $(1,0)$-vectors, one gets \begin{equation}\label{Schouten}
\{\oom\ot V, \oom'\ot V'\}=\oom'\we \iota_{V'}d\oom\ot V +\oom\we\iota_V
d\oom'\ot V'.\end{equation}

To construct deformations, we can now apply Kuranishi's recursive
formula like in \cite{MPPS}. We recall it for the sake of
completeness.

Let
$\{\be_1,\dots,\be_N\}$ be an orthonormal basis of the harmonic
representatives of $H^1(M,\Th_M)$. For any vector $\bt=(t_1,\dots, t_N)$
in $\bC^N$, let $\bmu(\bt)=t_1 \be_1+\cdots+t_N\be_N$ and set $\bph_1=\bmu$.
Then, one can define $\bph_r$ inductively for
$r\ge2$ as we will now recall.

Denote as usual by $\dbar^*$ the adjoint operator  to the
$\dbar$-operator on $M$ with respect to the Hermitian metric
previously defined and by $\triangle=\dbar\dbar^*+\dbar^*\dbar$
the Laplacian. Then we set \begin{equation}\label{phi}
\bph_r(\bt)=\half\sum_{s=1}^{r-1}\dbar^*\cg\{\bph_s(\bt),\bph_{rs}(\bt)\}=
\half\sum_{s=1}^{r-1}\cg\dbar^*\{\bph_s(\bt),\bph_{r-s}(\bt)\},\end{equation}
where $\cg$ is the corresponding Green's operator that inverts
$\triangle$ on the orthogonal complement of the space of harmonic
forms.

Consider the formal sum \begin{equation}\label{fsum}\bPh(\bt)=
\sum_{r\ge1}\bph_r.\end{equation} Observe that $\bPh$ belongs to
$\g^{*(0,1)}\ot \g^{1,0}$ and can be regarded either as a map
sending $\g^{0,1}$ to $\g^{1,0}$ or as one from $\g^{*(1,0)}$ to
$\g^{*(0,1)}$.

Let $\{\ga_1,\dots,\ga_P\}$ be an orthonormal basis for the space of harmonic
$(0, 2)$-forms with values in $\Th_M$. Define $f_k(\bt)$ to be the $L^2$-inner
product $\ll\{\bPh(\bt),\bPh(\bt)\},\ga_k\rr$. Kuranishi theory asserts the
existence of $\eps>0$ such that \begin{equation}\label{family}
\{\bt\in\bC^N:|\bt|<\eps,\ f_1(\bt)=0,\dots,f_P(\bt)=0\}\end{equation} forms a
locally complete family of deformations of $M$. We shall denote this set by
$\Kur$.

For each $\bt\in\Kur$, the associated sum $\bPh=\bPh(\bt)$ defines
a family of complex structures $J_{\bPh}$ whose (1,0)-forms are
given by $\omega-\bPh(\omega)$, $\omega \in \g^{*(1,0)}$ and whose
(0,1)-vectors are $\bar X + \bPh (\bar X)$, $\bar X \in \g^{0,1}$.

Consequently the integrability condition for a deformation $\bPh$
is, for any $\omega \in \g^{*(1,0)}$ and $\bar X, \bar Y\in
\g^{0,1}$,
\begin{equation}\label{integrability}
\left ( d(\omega-\bPh(\omega))\right )(\bar X+\bPh(\bar X), \bar
Y+\bPh(\bar Y))=0.
\end{equation}
To relate this condition to the Maurer-Cartan equation, one may
check that
\begin{equation}
-\left ( d(\omega-\bPh(\omega))\right )(\bar X+\bPh(\bar X), \bar
Y+\bPh(\bar Y))= \omega\left((\dbar \bPh+\frac12 \{\bPh, \bPh\})
(\bar X, \bar Y) \right).
\end{equation}


Now, given the fact that the harmonic theory is reduced to a
finite dimensional program as noted in Theorem \ref{harmonic rep},
a proof in \cite{MPPS} shows that every term in the power series
(\ref{fsum}) lies in $\g^{*(0,1)}\ot \g^{1,0}$. Thus

\begin{theorem}\label{main} Let $G$ be a  nilpotent Lie group with
co-compact subgroup $\Ga$, and let $J$ be an abelian invariant complex
structure on $M=\Ga\bs G$. Then the deformations arising from $J$
parameterized by (\ref{family}) are all invariant complex structures.
\end{theorem}

We can now find under which conditions $J_{\bPh}$ remains abelian.
Recall that a complex structure is abelian if and only if the
differential of a (1,0) form is of type (1,1). In other words, for
any $\omega \in \g^{*(1,0)}$ and $ X,  Y\in \g^{1,0}$,
\begin{equation}\label{abelian}
d(\omega-\bPh(\omega)) (X + \bPh (X), Y + \bPh (Y))=0.
\end{equation}
If one extends the Schouten-Nijenhuis bracket $\{\cdot ,\ \cdot\}$
to the exterior algebra by anti-derivative as seen in \cite{Poon}
and make use of the assumption that $J$ is abelian, a short
computation shows that
\begin{equation}\label{abelian condition}
d(\omega-\bPh(\omega)) (X + \bPh (X), Y + \bPh (Y))=\{ \bar{\bPh},
\omega-\bPh(\omega)\} (X,Y).
\end{equation}
Note that $\bPh(\omega)$ is in $\g^{*(0,1)}$, so $\{\bar{\bPh},
\bPh(\omega)\}$ is in $\g^{*(0,2)}$. Note further that if $\alpha$
is in $\g^{*(1,0)}$ and $\bar V$ is in $\g^{0,1}$, then
\begin{equation}
\{\alpha\otimes \bar V, \omega\}=\alpha\wedge\{\bar V,
\omega\}=\alpha\wedge \iota_{\bar V}d\omega.
\end{equation}
Since the complex structure is abelian, $d\omega$ is type (1,1).
Therefore, $\iota_{\bar V}d\omega$ is in $\g^{*(1,0)}$. Therefore,
$\{ \bar{\bPh}, \omega-\bPh(\omega)\}$ is in $\g^{*(2,0)}\otimes
\g^{1,0}$. It follows from equation (\ref{abelian condition})
above that the deformed complex structure $J_{\bPh}$ is abelian if
and only if
\begin{equation}\label{abelian condition 2}
\{ \bar{\bPh}, \omega-\bPh(\omega)\}=0
\end{equation}
for all $\omega\in \g^{*(1,0)}$.

\begin{theorem}\label{abeliandef}
$\bPh$ defines an abelian deformation if and only if it is
integrable and
$$\{ \bar{\bPh}, \omega-\bPh(\omega)\}=0 \qquad {\hbox{\rm for any }}  \omega \in \g^{*(1,0)}\,\, .$$
\end{theorem}

Infinitesimally, suppose $\bPh=t \bmu + t^2 \bph_2 + t^3 \bph_3 +
\dots$. Then looking at the degree one terms, the integrability
condition implies that $ \dbar \bmu = 0$. Computing the second
order term in the equation (\ref{abelian condition 2}) yields
\begin{equation}
\{\bmu, {\oomega}\}=0
\end{equation}
for any $\bar\omega$ in $\g^{*(0,1)}$. It follows that if $\bph_j$
are constructed through the Kuranishi recursive formula, then
$\bph_j=0$ for all $j\geq 2$. Conversely, if $\bmu$ is
$\dbar$-closed and satisfies the above condition, it is integrable
to an abelian complex structure. Therefore, we have the following
result.
\begin{proposition}
A parameter $\bmu\in H^1(M,\Th_M)$ defines an integrable
infinitesimal abelian deformation if and only if $\dbar \bmu = 0$
and
$$
\{ \bmu, \bar \omega \} =0 \qquad {\hbox{\rm for any }} \bar \omega \in \g^{*(0,1)}\, .
$$
\end{proposition}

Furthermore, we may consider $\bmu$ as an element in
$\g^{1,0}\oplus\g^{*(0,1)}$. With the Schouten bracket, the space
$\g^{1,0}\oplus\g^{*(0,1)}$ becomes a Lie algebra. Since the
complex structure is abelian, $\{\g^{1,0}, \g^{1,0}\}=0$. By
definition, $\{\g^{*(0,1)}, \g^{*(0,1)}\}=0$. Therefore, $\{ \bmu,
\bar \omega \} =0$  for any $\omega \in \g^{*(0,1)}$ if and only
if $\bmu$ is the in kernel of the adjoint map.

\begin{corollary} A parameter $\bmu\in H^1(M,\Th_M)$ defines an integrable
infinitesimal abelian deformation if and only if $\dbar \bmu = 0$
and $\bmu$ is in the center of the Lie algebra
$\g^{1,0}\oplus\g^{*(0,1)}$.
\end{corollary}

\section{Six dimensional examples} \label{6-examples}

We recall that a complex structure $J$ on a $2n$-dimensional
nilpotent Lie algebra $\mathfrak g$, one may consider
 the ascending series $\{ {\mathfrak g}_l^{J} \}$
defined inductively by ${\mathfrak g}_0^{J} = \{ 0 \}$ and
$$
{\mathfrak g}_l^{J} = \{ X \in {\mathfrak g}  \, : \, [X,
{\mathfrak g}] \subseteq {\mathfrak g}_{l-1}^{J}, \, [JX,
{\mathfrak g}] \subseteq {\mathfrak g}_{l-1}^{J} \}, \quad l \geq
1.
$$
The complex structure is said to be nilpotent if the series
satisfies ${\mathfrak g}_k^{J} = {\mathfrak g}$, for some positive
integer $k$ \cite{CFGU}. Apparently, an abelian complex structure
is nilpotent, with ${\mathfrak g}_\ell^{J} = {\mathfrak g}_\ell$,
for any $\ell \geq 0$.

By Corollary 2.7 in \cite{Ug} if $G$ is  a 6-dimensional nilpotent
Lie group admitting an invariant complex structure, then all
invariant complex structures are either nilpotent or non-nilpotent
altogether. Thus, combining this result with Theorem \ref{main}
one has the following observation.

\begin{corollary}\label{main2} Let $G$ be a  $6$-dimensional nilpotent
Lie group with co-compact subgroup $\Ga$, and let $J$ be an
abelian invariant complex structure on $M=\Ga\bs G$. Then the
deformations arising from $J$ parameterized by (\ref{family}) are
all invariant nilpotent complex structures.
\end{corollary}

Abelian complex structures on 2-step nilmanifolds with dimension
six are already studied in \cite{MPPS}.  Using the classification
obtained in \cite{simon,Ug}, in dimension 6 there are only two
$k$-step nilpotent Lie algebras with $k > 2$. Their respective
structure equations are as follow.
$$
\begin{array} {l}
{\fh}_9: [e_1, e_2] = e_3, \, [e_1, e_3] = [e_2, e_4] = e_6,\\
{\fh}_{15}: [e_1, e_2] = -e_4, [e_1, e_3] = [e_2, e_4] = e_5,
[e_1, e_4] = - [e_2, e_3] = - e_6.
\end{array}
$$
Both are $3$-step nilpotent.

It is known that any complex structure on the first one ${\fh}_9$
is abelian \cite[Theorem 2.9]{Ug}.  Thus, if $H_9$ is the simply
connected nilpotent Lie group with Lie algebra ${\fh}_9$ and
  $\Ga$ is a co-compact subgroup $\Ga$,  the deformations arising from
any invariant complex structure $J$ on $\Ga \backslash
H_9$ parameterized by (\ref{family}) are all invariant abelian
complex structures.

We can actually verify directly that, if we start from the abelian complex structure $J$ on ${\fh}_9$ with
$$
Je_1=e_2\, , \qquad Je_3 =e_4\, , \qquad Je_5=e_6\, ,
$$
any Kuranishi deformation is still abelian. More explicitly,
$$
X_1= e_5-ie_6\, , \quad
X_2= e_3-ie_4\, , \quad
X_3= e_1-ie_2\,
$$
is a basis of left invariant (1,0) fields.
Moreover
$$
\bar \omega^1= e^5-ie^6\, , \quad
\bar \omega^2= e^3-ie^4\, , \quad
\bar \omega^3= e^1-ie^2\, ,
$$
where $e^i$ are the dual of $e_j$ gives a basis of $(0,1)$-forms. An orthonormal  basis of
harmonic representatives for $H^1_{\dbar} (\fh_9^{1,0})$ is given by
\[
\beta_1= \bar \omega^1 \ot X_1\, , \quad
\beta_2= \frac 1 {\sqrt 2} (\bar \omega^2 \ot X_2 -  \bar \omega^3 \ot X_3)\, , \quad
\beta_3= \frac 1 {\sqrt 2} (\bar \omega^2 \ot X_1 -  \bar \omega^3 \ot X_2)\, .
\]
Now if $ \bmu=a_1\beta_1+ a_2\beta_2+a_3\beta_3$ one can check
that $\{\bmu, \bar \omega^i\}=0$ for any $i=1,2,3$, which shows
directly that any small Kuranishi deformation of $\fh_9$ arising
from $J$ is abelian.

We also verify directly that Kuranishi deformations of abelian complex structures on $\fh_{15}$ are not necessarily abelian.
Let us start from the abelian complex structure $J$ on ${\fh}_{15}$ with
$$
Je_1=e_2\, , \qquad Je_3 =e_4\, , \qquad Je_5=e_6\, .
$$
Like in the previous case,
\begin{eqnarray*}
&X_1= e_5-ie_6\, , \quad
X_2= e_3-ie_4\, , \quad
X_3= e_1-ie_2\, \\
& \bar \omega^1= e^5-ie^6\, , \quad
\bar \omega^2= e^3-ie^4\, , \quad
\bar \omega^3= e^1-ie^2\, ,
\end{eqnarray*}
are basis of left invariant (1,0) fields and $(0,1)$-forms respectively, where $e^i$ are the dual of $e_j$. An orthonormal  basis of
harmonic representatives for $H^1_{\dbar} (\fh_{15}^{1,0})$ is given by
\begin{eqnarray*}
\beta_1&=&  \oomega^1 \ot X_1\, , \quad
\beta_2= \frac 1 {\sqrt 5} (\oomega^2 \ot X_1 -  2 \oomega^3 \ot X_2)\, , \quad
\beta_3=\oomega^3 \ot X_1 \, , \\
\beta_4&=&\oomega^1 \ot X_2 \, , \quad
\beta_5=\oomega^2 \ot X_2\, .
\end{eqnarray*}
Setting $\bmu=\sum_{\ell=1}^5 a_\ell \beta_\ell$, we have $\{\bmu,
\bar \omega^\ell\}=0$ for any $\ell=1,2,3$ if and only if $a_4=a_5
= 0$. Hence any such $\bmu$ with non vanishing $a_4$ or $a_5$
parametrizes deformations of $J$ which are invariant and
nilpotent, but not abelian.

\section{Higher Dimension Example}\label{instability}

In higher dimensions, deformations of abelian complex structures are not necessarily nilpotent.
Here is an  example of an abelian complex structure on a ten dimensional 3-step nilpotent Lie algebra $\fn$ which deforms into non-nilpotent complex structures.

The nilpotent Lie algebra $\fn$ has structure equations:
\begin{eqnarray*}
& &d e^1 = d e^2 = d e^3 = d e^7 = 0, \quad d e^4 =-e^{12}+e^{13}+e^{27}, \quad   d e^5 = -e^{12}-e^{17}+e^{23}, \\
 & & d e^6 =- e^{14} - e^{15} - e^{25} +e^{24} -e^{19} +e^{28}-2e^{45}+e^{48}+e^{59}-e^{49}+e^{58}-e^{89}, \\
& &d e^8 = -e^{17}+e^{23}-e^{13}-e^{27}, \quad d e^9 =
2e^{12}+e^{17}-e^{23}-e^{13}-e^{27}, \quad d e^{10} = 0.
\end{eqnarray*}

Considers the family of complex structure $J_{s,t}$ on $\fn$ such
that
\begin{eqnarray*}
& &J_{s,t} e_1 = e_2, J_{s,t} e_4 = e_5, J_{s,t} e_8 = e_9,\\
& &J_{s, t} e_3 = t e_6 + s e_7, J_{s,t} e_{10} = - s e_6 - t
e_7,\\
& & J_{s, t} e_6 = \frac {1}{t^2 -s^2} (- t e_3 - s e_{10}), J_{s, t}
e_7 = \frac {1}{t^2 -s^2}(s e_3 + t e_{10})
\end{eqnarray*}
 with $s,t$ real parameters such that $t^2
\neq s^2$. To check integrability, choose the (1,0)-forms as
follows.
\begin{eqnarray*}
\omega^1 &=& e^1+ie^2, \quad \omega^2= e^4+ie^5, \quad \omega^3= e^8+ie^9, \\
\omega^4&=& e^3+i(te^6+se^7), \quad \omega^5= e^{10}-i(se^6+te^7).
\end{eqnarray*}
It suffices to show that $d\omega^j$ is of type (1,1)+(2,0), i.e.,
it is generated by the ideal of type (1,0)-forms. It is a long,
yet straight forward exercise. We do not display all formula here.

For  $t = 0, s = 1$ one has an abelian complex structure. Note
that the center of $\fn$ is spanned by $e_6$ and $e_{10}$. If a
complex structure $J$ is nilpotent, then $\g_1^J$ is non-trivial,
$J$-invariant and contained in the center of $\g$. Given the
dimension of the center, it is possible only if the center of $\g$
is equal to $\g_1^J$. In particular, it is $J$-invariant. Note
that the center is not preserved by $J_{s,t}$ if $t \neq 0$. Thus
a generic $J_{s,t}$ is not nilpotent.

%

Dipartimento di Matematica, Universit\`a di Torino, via Carlo Alberto 10,
I-10128 Torino, Italy (\texttt{sergio.console@unito.it})

Dipartimento di Matematica, Universit\`a di Torino, via Carlo
Alberto 10, I-10128 Torino, Italy
(\texttt{annamaria.fino@unito.it})

Department of Mathematics, University of California at Riverside, Riverside,
CA 92521, USA (\texttt{ypoon@math.ucr.edu})

\enddocument